\documentstyle[fleqn,12pt]{article}
\begin{document}
\pagenumbering{arabic}
\setcounter{page}{1}
\pagestyle{plain}
\baselineskip=18pt

\thispagestyle{empty}
\rightline{YTUMB 2001-03, November 2001} 
\vspace{1.4cm}

\begin{center}
{\Large\bf Z$_3$-graded differential geometry of \\ quantum plane}
\end{center}

\vspace{1cm}
\begin{center} Salih Celik \footnote{E-mail: sacelik@yildiz.edu.tr}
 
Yildiz Technical University, Department of Mathematics, \\
34210 DAVUTPASA-Esenler, Istanbul, TURKEY. \end{center}

\vspace{2cm}
{\bf Abstract}

In this work, the Z$_3$-graded differential geometry of the quantum plane 
is constructed. The corresponding quantum Lie algebra and its Hopf 
algebra structure are obtained. The dual algebra, i.e. universal enveloping algebra of the quantum plane is explicitly constructed and an isomorphism between the quantum Lie algebra and the dual algebra is given. 

\vfill\eject
\noindent
{\bf 1. Introduction}

\noindent
After the discovery of the quantum plane by Manin [1], Wess and Zumino [2] 
developed a differential calculus on the quantum (hyper)plane covariant with 
respect to the action of the quantum group. In their method, the R-matrix is 
obtained using the consistency conditions. This leads to a consistent exterior 
derivative. The purely algebraic properties of these recently discovered 
spaces have been deeply discussed. The $q$-differential algebras have become 
the object of excellent works [3,4]. 

The Z$_3$-graded algebraic structures have been introduced and studied [5]. 
The de Rham complex with differential operator {\sf d} satisfying the 
$Q$-Leibniz rule and the condition ${\sf d}^3 = 0$ on an associative unital 
algebra has been constructed by Bazunova et al [6] using the methods of 
Ref. 2. This paper consider an alternative approach where, instead of adopting 
an R-matrix, consistency conditions on natural commutation relations are used. 

The cyclic  group Z$_3$ can be represented in the complex plane by means of 
the cubic roots of 1: let $j = e^{{2\pi i}\over 3}$ $(i^2 = - 1)$. Then one 
has 
$$j^3 = 1 \quad \mbox{and} \quad j^2 + j +1 = 0 \quad \mbox{or} \quad 
  (j + 1)^2 = j.  \eqno(1)$$
One can define the Z$_3$-graded commutator $[A,B]$ as [3]
$$[A,B]_{Z_3} = AB - j^{deg(A) deg(B)} BA \eqno(2)$$
where $deg(X)$ denotes the grade of $X$. If $A$ and $B$ are $j$-commutative, 
then we have 
$$AB = j^{deg(A) deg(B)} BA. \eqno(3)$$

\vspace*{0.3cm}
\noindent
{\bf 2. Review of Hopf algebra ${\cal A}$} 

\noindent
Elementary properties of the extended quantum plane are described in Ref. 7. 
We state briefly the properties we are going to need in this work. 

\noindent
{\bf 2.1 The algebra of polynomials on the $q$-plane}

\noindent
The quantum plane [1] is defined as an associative algebra generated by two 
noncommuting coordinates $x$ and $y$ with the relation 
$$ x y - q y x = 0 \qquad q \in {\cal C} - \{0\}. \eqno(4) $$
This associative algebra over the complex numbers ${\cal C}$, is known 
as the algebra of polynomials over the quantum plane and is often denoted by 
$C_q[x,y]$. In the limit $q \longrightarrow 1$, this algebra is commutative 
and can be considered as the algebra of polynomials $C[x,y]$ over the usual 
plane, where $x$ and $y$ are the two coordinate functions. In below, we show 
that there cannot exits a Z$_3$-graded {\it commutative} differential calculus 
as in the Z$_2$-grade case. We denote the unital extension of $C_q$ by 
${\cal A}$. 

\noindent
{\bf 2.2 Hopf algebra structure on ${\cal A}$}

\noindent
The definitions of a coproduct, a counit and a coinverse on the algebra 
${\cal A}$ are as follows [7,8]: 

{\bf (1)} The coproduct 
$\Delta_{\cal A} : {\cal A} \longrightarrow {\cal A} \otimes {\cal A}$ 
is defined by 
$$\Delta_{\cal A}(x) = x \otimes x \qquad 
  \Delta_{\cal A}(y) = y \otimes 1 + x \otimes y \eqno(5)$$
is coassociative: 
$$(\Delta_{\cal A} \otimes \mbox{id}) \circ \Delta_{\cal A} = 
  (\mbox{id} \otimes \Delta_{\cal A}) \circ \Delta_{\cal A} \eqno(6)$$
where id denotes the identity map on ${\cal A}$. 

{\bf (2)} The counit $\epsilon_{\cal A} : {\cal A} \longrightarrow {\cal C}$ 
is given by 
$$\epsilon_{\cal A}(x) = 1 \qquad \epsilon_{\cal A}(y) = 0. \eqno(7)$$
The counit $\epsilon_{\cal A}$ has the property 
$$m_{\cal A} \circ (\epsilon_{\cal A} \otimes \mbox{id}) \circ \Delta_{\cal A} 
= m_{\cal A} \circ (\mbox{id} \otimes \epsilon_{\cal A}) \circ \Delta_{\cal A} 
\eqno(8)$$
where $m_{\cal A}$ stands for the algebra product 
${\cal A} \otimes {\cal A} \longrightarrow {\cal A}$. 

{\bf (3)} If we extend the algebra $\cal A$ by adding the inverse of $x$ 
then the algebra $\cal A$ admits a ${\cal C}$-algebra antihomomorphism 
(coinverse) 
$\kappa_{\cal A} : {\cal A} \longrightarrow {\cal A}$ defined by 
$$\kappa_{\cal A}(x) = x^{-1} \qquad \kappa_{\cal A}(y) = - x^{-1} y.\eqno(9)$$
The coinverse $\kappa_{\cal A}$ satisfies 
$$m_{\cal A} \circ (\kappa_{\cal A} \otimes \mbox{id}) \circ \Delta_{\cal A} = 
  \epsilon_{\cal A} = 
  m_{\cal A} \circ (\mbox{id} \otimes \kappa_{\cal A}) \circ \Delta_{\cal A}. 
\eqno(10)$$

\noindent
{\bf 3. Construction of bicovariant Z$_3$-graded differential calculus 
on $\cal A$} 

\noindent
The Woronowicz theory [9] is based on the idea that the differential and 
algebraic structures of $\cal A$ can coact covariantly on the algebra of its 
differential calculus over $\cal A$. We first recall some basic notions about 
differential calculus on the extended $q$-plane. 

\noindent
{\bf 3.1 Differential algebra} 

\noindent
To begin with, we note the properties of the exterior differential {\sf d}. 
The exterior differential {\sf d} is an operator which gives the mapping from 
the generators of ${\cal A}$ to the differentials: 
$${\sf d} : a \longrightarrow {\sf d} a \qquad a \in \{x,y\}. $$
We demand that the exterior differential {\sf d} has to satisfy two 
properties: 
$${\sf d}^3 = 0 \eqno(11)$$
and the Z$_3$-graded Leibniz rule 
$${\sf d}(f g) = ({\sf d} f) g + j^{deg(f)} ({\sf d} g). \eqno(12)$$

In order to establish a noncommutative differential calculus including second 
order differentials of the generators of $\cal A$ on the $q$-plane, we assume 
that the commutation relations between the coordinates and their first order 
differentials are of the following form: 
$$x ~{\sf d}x = A {\sf d}x ~x $$
$$x ~{\sf d}y = C_{11} {\sf d}y ~x + C_{12} {\sf d}x ~y $$
$$y ~{\sf d}x = C_{21} {\sf d}x ~y + C_{22} {\sf d}y ~x $$
$$y ~{\sf d}y = B {\sf d}y ~y. \eqno(13)$$

The coefficients $A$, $B$ and $C_{ik}$ will be determined in terms of the 
complex deformation parameter $q$ and $j$. To find them we shall use the 
covariance of the noncommutative differential calculus. 

Since we assume that ${\sf d}^3 = 0$ and ${\sf d}^2 \neq 0$, in order to 
construct a self-consistent theory of differential forms it is necessary to 
add to the first order differentials of coordinates ${\sf d} x$, 
${\sf d} y$ a set of {\it second order differentials} ${\sf d}^2 x$, 
${\sf d}^2 y$. Let us begin by assuming that 
$${\sf d} x ~{\sf d} y = F {\sf d} y ~{\sf d} x \qquad 
  ({\sf d}x)^3 = 0 = ({\sf d}y)^3 \eqno(14)$$
where $F$ is a parameter that shall described later. 

The first differentiation of (13) gives rise to the relations between the 
generators $x$, $y$ and second order differentials ${\sf d}^2 x$,  
${\sf d}^2 y$ including first order differentials: 
$$x ~{\sf d}^2 x = A ~{\sf d}^2 x ~x + (Aj - 1) ({\sf d} x)^2 $$
$$x ~{\sf d}^2 y = C_{11} {\sf d}^2 y ~x + C_{12} {\sf d}^2 x ~y 
   + K_1 {\sf d} y ~{\sf d} x $$
$$y ~{\sf d}^2 x = C_{21} {\sf d}^2 x ~y + C_{22} {\sf d}^2 y ~x + 
  K_2 {\sf d}y ~{\sf d} x $$ 
$$y ~{\sf d}^2 y = B {\sf d}^2 y ~y + (Bj - 1) ({\sf d} y)^2. \eqno(15) $$
where 
$$K_1 = j C_{11} +j C_{12} F - F \qquad K_2 = j C_{21} F + j C_{22} - 1. 
\eqno(16)$$
The relations (15) are {\it not homogeneous} in the sense that the 
commutation relations between the generators and second order differentials 
include first order differentials as well. In the following subsection, we 
shall see that the commutation relations between the coordinates and their 
second order differentials can be made homogeneous. They will not include 
first order differentials by removing them using the covariance of the 
noncommutative differential calculus. 

Applying the exterior differential {\sf d} to the relations (15), we get 
$${\sf d} x ~{\sf d}^2 x = j^{-2} {\sf d}^2 x ~{\sf d} x $$
$${\sf d}x ~{\sf d}^2 y = j^2 C_{11} Q_1^{-1} {\sf d}^2 y ~{\sf d}x + 
  (j^2 C_{12} + F^{-1} K_1) Q_1^{-1} {\sf d}^2 x ~{\sf d} y $$
$${\sf d} y ~{\sf d}^2 x = j^2 C_{21} Q_2^{-1} {\sf d}^2 x ~{\sf d} y 
   + (j^2 C_{22} + K_2) Q_2^{-1} {\sf d}^2 y ~{\sf d} x $$ 
$${\sf d} y ~{\sf d}^2 y = j^{-2} {\sf d}^2 y ~{\sf d} y \eqno(17)$$
where 
$$Q_1 = - j^2 (C_{12} + C_{11} F^{-1} + 1) \qquad 
  Q_2 = - j^2 (C_{22} + C_{21} F + 1). \eqno(18)$$
The differentiation of second or third relations of (17) gives rise to the 
relations between the second order differentials: 
$${\sf d}^2 x ~{\sf d}^2 y = F {\sf d}^2 y {\sf d}^2 x.  \eqno(19)$$

\vspace*{0.3cm}\noindent
{\bf 3. 2 Covariance}

\noindent
In order to homogenize the relations (15), we shall consider the covariance of 
the noncommutative differential calculus. Let $\Gamma$ be a bimodule over the 
algebra $\cal A$ generated by the elements of the set 
$\{x, y, {\sf d} x, {\sf d} y, {\sf d}^2 x, {\sf d}^2 y\}$. One says that 
$(\Gamma, {\sf d})$ is a {\it first-order differential calculus} over the 
Hopf algebra $({\cal A}, \Delta_{\cal A},\epsilon_{\cal A}, \kappa_{\cal A})$. 
We begin with the definitions of a left- and right-covariant bimodule. 

{\bf (1)} Let $\Gamma$ be a bimodule over ${\cal A}$ and 
$\Delta^R : \Gamma \longrightarrow \Gamma \otimes {\cal A}$ 
be a linear homomorphism. We say that $(\Gamma, \Delta^R)$ is a 
right-covariant bimodule if 
$$\Delta^R(a \rho + \rho' a') = \Delta_{\cal A}(a) \Delta^R(\rho) + 
  \Delta^R(\rho') \Delta_{\cal A}(a') \eqno(20)$$
for all $a, a' \in {\cal A}$ and $\rho, \rho' \in \Gamma$, and 
$$(\Delta^R \otimes \mbox{id}) \circ \Delta^R = 
  (\mbox{id} \otimes \Delta_{\cal A}) \circ \Delta^R \qquad 
  (\mbox{id} \otimes \epsilon) \circ \Delta^R = \mbox{id}. \eqno(21)$$
The action of $\Delta^R$ on the first order differentials is 
$$\Delta^R({\sf d}x) = {\sf d}x \otimes x \qquad 
  \Delta^R({\sf d}y) = {\sf d}y \otimes 1 + {\sf d}x \otimes y \eqno(22)$$
since 
$$\Delta^R({\sf d}a) = ({\sf d} \otimes \mbox{id}) \Delta_{\cal A}(a) \quad 
  \forall a \in {\cal A}. \eqno(23)$$

We now apply the linear map $\Delta^R$ to relations (13): 
$$\Delta^R(x ~{\sf d} x) = \Delta_{\cal A}(x) \Delta^R ({\sf d} x) = 
  A \Delta^R ({\sf d} x ~x), $$
$$\Delta^R(x ~{\sf d} y) = C_{11} \Delta^R ({\sf d} y ~x) + 
   C_{12} \Delta^R({\sf d} x ~y) + 
   (qA - C_{11} - qC_{12}) {\sf d}x ~x \otimes x y, $$
$$\Delta^R(y ~{\sf d} x) = C_{21} \Delta^R ({\sf d} x ~y) + 
   C_{22} \Delta^R({\sf d} y ~x) + 
   (A - q C_{21} - C_{22}) {\sf d}x ~x \otimes y x, $$
\begin{eqnarray*}
\Delta^R(y ~{\sf d} y) 
& = & B \Delta^R ({\sf d} y ~y) + 
     (C_{12} + C_{21} - B) {\sf d}x ~y \otimes y + 
     (A - B) {\sf d} x ~x \otimes y^2 \\
& & + (C_{11} + C_{22} - B) {\sf d} y ~x \otimes y, 
\end{eqnarray*}
and to relations (14) 
$$\Delta^R({\sf d} x ~{\sf d} y) = F \Delta^R ({\sf d} y ~{\sf d} x) + 
   (q - F) ({\sf d} x)^2 \otimes y x. $$
So we must have 
$$C_{11} + q C_{12} = q A \qquad C_{11} + C_{22} = B \qquad A = B $$
$$q C_{21} + C_{22} = A \qquad C_{12} + C_{21} = B \qquad F = q. \eqno(24)$$

{\bf (2)} Let $\Gamma$ be a bimodule over ${\cal A}$ and 
$\Delta^L : \Gamma \longrightarrow {\cal A} \otimes \Gamma$ 
be a linear homomorphism. We say that $(\Gamma, \Delta^L)$ is a 
left-covariant bimodule if 
$$\Delta^L(a \rho + \rho' a') = \Delta_{\cal A}(a) \Delta^L(\rho) + 
  \Delta^L(\rho') \Delta_{\cal A}(a') \eqno(25)$$
for all $a, a' \in {\cal A}$ and $\rho, \rho' \in \Gamma$, and 
$$(\Delta_{\cal A} \otimes \mbox{id}) \circ \Delta^L = 
  (\mbox{id} \otimes \Delta^L) \circ \Delta^L \qquad 
  (\epsilon \otimes \mbox{id}) \circ \Delta^L = \mbox{id}. \eqno(26)$$
Since 
$$\Delta^L({\sf d}a) = (\mbox{id} \otimes {\sf d}) \Delta_{\cal A}(a) \quad 
  \forall a \in {\cal A} \eqno(27)$$
the action of $\Delta^L$ on the first order differentials gives rise to the 
relations 
$$\Delta^L({\sf d}x) = x \otimes {\sf d} x \qquad 
  \Delta^L({\sf d}y) = x \otimes {\sf d} y. \eqno(28)$$
Applying $\Delta^L$ to relations (13), we get 
$$C_{12} = 0 \qquad C_{21} = q^{-1} \qquad B = q^{-1}. \eqno(29)$$
With the relations (24), we then obtain 
$$A = q^{-1} \qquad C_{11} = 1 \qquad C_{21} = q^{-1} $$
$$B = q^{-1} \qquad C_{12} = 0 \qquad C_{22} = q^{-1} - 1. \eqno(30)$$
So 
$$K_1 = j - q \qquad K_2 = q^{-1}(j - q) \qquad 
  Q_1 = - j^2 (q^{-1} + 1) = Q_2.$$

On the other hand, since the differential of a function $f$ of the coordinates 
$x$ and $y$ is of the form 
$${\sf d} f = ({\sf d} x \partial_x + {\sf d} y \partial_y) f \eqno(31)$$
and 
$${\sf d}^2 f = \left({\sf d}^2x \partial_x + {\sf d}^2y \partial_y + 
  j ({\sf d}x)^2 \partial_x^2 + j ({\sf d}y)^2 \partial_y^2 + 
 {\sf d}x {\sf d}y (\partial_x \partial_y + q \partial_y \partial_x)\right)f,$$
\begin{eqnarray*}
{\sf d}^3 f 
& = & {\sf d}^2x {\sf d}y (j^2 \partial_y \partial_x + 
      q^{-1} j \partial_x \partial_y + 
      {{1-jq^{-1}}\over {q+1}} \partial_x \partial_y + 
      {{q-j}\over {q+1}} \partial_y \partial_x) f \\
&   & + {\sf d}^2y {\sf d}x (j^2 \partial_x \partial_y - 
       {{j^2}\over {q+1}} \partial_x \partial_y - 
      {1 \over {q+1}} \partial_y \partial_x) f + \cdots \\
& = & {{j^2}\over {q+1}} {\sf d}^2x {\sf d}y (\partial_y \partial_x - 
      \partial_x \partial_y) + {1 \over {q+1}} {\sf d}^2y {\sf d}x 
      (q j^2 \partial_x \partial_y - \partial_y \partial_x) f + \cdots \\
&\equiv & 0
\end{eqnarray*}
we have 
$$\partial_x \partial_y = \partial_y \partial_x \eqno(32)$$
if $q$ satisfies the identities 
$$q j^2 = 1 \qquad q^2 + q + 1 = 0. \eqno(33)$$
One can then chose 
$$q = j^{-2} = j. \eqno(34)$$
Consequently, the relations (13)-(15), (17) and (19) are explicity as follows: 
the commutation relations of the coordinates and their first order 
differentials are [10] 
$$x~ {\sf d}x = q^{-1} ~{\sf d}x~ x \qquad 
  x~ {\sf d}y = {\sf d} y~ x $$
$$y~ {\sf d}x = q^{-1} {\sf d}x~ y + (q^{-1} - 1) {\sf d}y ~x \qquad 
  y~ {\sf d}y = q^{-1} {\sf d}y~ y \eqno(35)$$
and among those first order differentials are 
$${\sf d} x ~{\sf d} y = q ~{\sf d} y ~{\sf d} x 
  \qquad ({\sf d} x)^3 = 0 = ({\sf d} y)^3. \eqno(36)$$

The commutation relations between variables and second order differentials are 
$$x ~{\sf d}^2 x = q^{-1} {\sf d}^2 x ~x \qquad 
  x ~{\sf d}^2 y = {\sf d}^2 y ~x $$
$$y ~{\sf d}^2 y = q^{-1} {\sf d}^2 y ~y \qquad 
  y ~{\sf d}^2 x = q^{-1} {\sf d}^2 x ~y + (q^{-1} - 1) {\sf d}^2 y ~x. 
  \eqno(37) $$

The commutation relations between first order and second order differentials 
are 
$${\sf d} x ~{\sf d}^2 x = q^{-2} {\sf d}^2 x ~{\sf d} x \qquad 
  {\sf d}x ~{\sf d}^2 y = q^2 {\sf d}^2 y ~{\sf d} x $$
$${\sf d} y ~{\sf d}^2 y = q^{-2} {\sf d}^2 y ~{\sf d} y \qquad 
  {\sf d} y ~{\sf d}^2 x = q^{-2} {\sf d}^2 x ~{\sf d} y + 
  (q - q^{-1}) ~{\sf d}^2 y ~{\sf d} x \eqno(38)$$
and those among the second order differentials are 
$${\sf d}^2 x ~{\sf d}^2 y = q~ {\sf d}^2 y ~{\sf d}^2 x.   \eqno(39)$$

Now, it can be checked that the linear maps $\Delta^R$ and $\Delta^L$ leave 
invariant the relations (35)-(39). One can also check that the identities 
(21), (26) and also the following identities are satisfied: 
$$(\mbox{id} \otimes {\sf d}) \Delta_{\cal A}(a) = \Delta^L({\sf d} a) \qquad 
 ({\sf d} \otimes \mbox{id}) \Delta_{\cal A}(a) = \Delta^R({\sf d} a) 
\eqno(40)$$
and 
$$(\Delta^L \otimes \mbox{id}) \circ \Delta^R = 
  (\mbox{id} \otimes \Delta^R) \circ \Delta^L. \eqno(41)$$

\vspace*{0.3cm}\noindent
{\bf 4. Cartan-Maurer one-forms on ${\cal A}$}

\noindent
In analogy with the left-invariant one-forms on a Lie group in classical 
differential geometry, one can construct two one-forms using the generators of 
$\cal A$ as follows [7]: 
$$\theta = {\sf d}x ~x^{-1} \qquad 
  \varphi = {\sf d} y - {\sf d} x ~x^{-1} y. \eqno(42)$$
The commutation relations between the generators of $\cal A$ and one-forms 
are [7] 
$$x \theta = q^{-1} \theta x \qquad 
  y \theta = q^{-1} \theta y + (q^{-1} - 1) \varphi $$
$$x \varphi = \varphi x \qquad y \varphi = \varphi y. \eqno(43)$$
The first order differentials with one-forms satisfy the following relations 
$$\theta ~{\sf d} x = q {\sf d} x ~\theta \qquad 
  \varphi ~{\sf d} x = {\sf d} x ~\varphi $$
$$\theta ~{\sf d} y = q {\sf d} y \theta \qquad 
  \varphi ~{\sf d} y = {\sf d} y ~\varphi \eqno(44)$$
and with second order differentials 
$$\theta ~{\sf d}^2 x = q^2 {\sf d}^2 x ~\theta \qquad 
  \theta ~{\sf d}^2 x = q^2 {\sf d}^2 x ~\theta $$
$$\varphi ~{\sf d}^2 x = q^{-2} {\sf d}^2 x ~\varphi \qquad 
  \varphi ~{\sf d}^2 y = q^{-2} {\sf d}^2 y ~\varphi. \eqno(45)$$

The commutation rules of the elements $\theta$ and $\varphi$ are 
$$\theta^3 = 0 \qquad \theta \varphi = \varphi \theta \eqno(46\mbox{a})$$
and
$$\varphi^3 = 0 \eqno(46\mbox{b})$$
provided that $q^2 + q + 1 = 0.$ 

We denote the algebra of the forms generated by the two elements $\theta$ and 
$\varphi$ by $\Omega$. We make the algebra $\Omega$ into a Z$_3$-graded Hopf 
algebra with the following co-structures [7]: 
the coproduct 
$\Delta_\Omega : \Omega \longrightarrow \Omega \otimes \Omega$ is defined by 
$$\Delta_\Omega(\theta) = \theta \otimes 1 + 1 \otimes \theta \qquad 
  \Delta_\Omega(\varphi) = 
   \varphi \otimes 1 + x \otimes \varphi - y \otimes \theta.   \eqno(47) $$
The counit $\epsilon_\Omega : \Omega \longrightarrow {\cal C}$ is given by 
$$\epsilon_\Omega(\theta) = 0 \qquad \epsilon_\Omega(\varphi) = 0 \eqno(48)$$
and the coinverse $\kappa_\Omega: \Omega \longrightarrow \Omega$ is defined by 
$$\kappa_\Omega(\theta) = - \theta \qquad 
\kappa_\Omega(\varphi) = - q^{-1} \varphi x^{-1} - \theta x^{-1} y. \eqno(49)$$
One can easily check that (6), (8) and (10) are satisfied. Note that 
the commutation relations (43)-(46) are compatible with $\Delta_\Omega$, 
$\epsilon_\Omega$ and $\kappa_\Omega$, in the sense that 
$\Delta_\Omega(x \theta) = q^{-1} \Delta_\Omega(\theta x)$, and so on. 

\vfill\eject
\noindent
{\bf 5. Quantum Lie algebra} 

\noindent
The commutation relations of Cartan-Maurer forms allow us to construct the 
algebra of the generators. In order to obtain the quantum Lie algebra 
of the algebra generators we first write the Cartan-Maurer forms as 
$${\sf d} x = \theta x \qquad {\sf d} y = \varphi + \theta y. \eqno(50)$$
The differential {\sf d} can then the expressed in the form 
$${\sf d} = \theta H + \varphi X. \eqno(51)$$
Here $H$ and $X$ are the quantum Lie algebra generators. 
We now shall obtain the commutation relations of these generators. 
Considering an arbitrary function $f$ of the coordinates of the quantum plane 
and using that ${\sf d}^3 = 0$ one has 
$${\sf d}^2 f = {\sf d} \theta ~H f + {\sf d} \varphi ~X f + 
   j \theta ~{\sf d} H f + j \varphi ~{\sf d} X f, $$
and 
$${\sf d}^3 f = {\sf d}^2 \theta ~H f + {\sf d}^2 \varphi ~X f + 
             j^2 {\sf d} \theta ~{\sf d} H f + j^2 {\sf d} \varphi ~{\sf d} Xf 
             + j^2 \theta ~{\sf d}^2 ~H f + j^2 \varphi ~{\sf d}^2 X f.$$

So we need the two-forms. Applying the exterior differential {\sf d} to the 
relations (42) one has 
$${\sf d} \theta = {\sf d}^2 x ~x^{-1} - j \theta^2 $$
$${\sf d} \varphi = {\sf d}^2 y - {\sf d}^2 x ~x^{-1} y - j \theta \varphi. 
\eqno(52)$$
Also, since 
$$\theta ~{\sf d} \theta = q^{-2} {\sf d} \theta ~\theta $$
$$\theta ~{\sf d} \varphi = q^2 {\sf d} \varphi ~\theta + 
        (q - q^{-1}) {\sf d} \theta ~\varphi + (q^{-1} - q) \theta^2 \varphi $$
$$\varphi ~{\sf d} \theta = q^{-2} {\sf d} \theta ~\varphi + 
                             (q^{-1} - q) \theta^2 \varphi $$
$$\varphi ~{\sf d} \varphi = q^{-2} {\sf d} \varphi ~\varphi + 
  (q^{-1} - q) \theta \varphi^2 \eqno(53)$$
we have 
$${\sf d}^2 \theta = 0 \qquad 
 {\sf d}^2 \varphi = j {\sf d} \theta ~\varphi - j {\sf d}\varphi ~\theta - 
 j \theta^2 \varphi. \eqno(54)$$
Using the Cartan-Maurer equations we find the following commutation 
relations for the quantum Lie algebra: 
$$X H = q^{-1} H X + X. \eqno(55)$$

The commutation relation (55) of the algebra generators should be consistent 
with monomials of the coordinates of the quantum plane. To do this, we 
evaluate the commutation relations between the generators of 
algebra and the coordinates. The commuation relations of the generators with 
the coordinates can be extracted from the Z$_3$-graded Leibniz rule: 
\begin{eqnarray*}
{\sf d} (x f) 
& = & ({\sf d} x) f + x ({\sf d} f) \\
& = & \theta (x + q^{-1} x H) f + \varphi (x X) f \\ 
& = & (\theta H + \varphi X) x f \hspace*{7.8cm}{(56)}
\end{eqnarray*}
and 
\begin{eqnarray*}
{\sf d} (y f) 
& = & ({\sf d} y) f + y ({\sf d} f) \\
& = & \theta (y + q^{-1} y H) f + \varphi (1 + y X + (q^{-1} - 1) H) f \\ 
& = & (\theta H + \varphi X) y f. \hspace*{7.7cm}{(57)}
\end{eqnarray*}
This yields 
$$H x  = x + q^{-1} x H \qquad H y = y + q^{-1} y H $$
$$X x = x X \qquad X y = 1 + yX + (q^{-1} - 1) H.   \eqno(58)$$

We know that the differential operator {\sf d} satisfies the Z$_3$-graded 
Leibniz rule. Therefore, the generators $H$ and $X$ are endowed with a 
natural coproduct. To find them, we need to the following commutation 
relation 
$$H x^m = {{1 - q^{-m}}\over {1-q^{-1}}} x^m + q^{-m} x^m H \eqno(59\mbox{a})$$
and 
$$H y^n = {{1- q^{-n}}\over {1 - q^{-1}}} y^n + q^{-n} y^n H\eqno(59\mbox{b})$$
where use was made of (58). The relation (59a) is understood as an operator 
equation. This implies that when $H$ acts on arbitrary monomials $x^m y^n$, 
$$H (x^m y^n) = {{1 - q^{-(m + n)}}\over {1 - q^{-1}}} (x^m y^n) + 
  q^{-(m+n)} (x^m y^n) H \eqno(60)$$
from which we obtain 
$$H = { {1 - q^{-N}} \over {1 - q^{-1}}}  \eqno(61)$$
where $N$ is a number operator acting on a monomial as 
$$N(x^m y^n) = (m + n) x^m y^n. \eqno(62)$$
We also have 
$$X(x^m y^n) = (x^m y^n) X + {{1 - q^{-n}}\over {1 - q^{-1}}} x^m y^{n-1} 
  (1 + (q^{-1} - 1) H). \eqno(63)$$

So, applying the Z$_3$-graded Leibniz rule to the product of functions $f$ and 
$g$, we write 
$${\sf d} (f g) = [(\theta H + \varphi X) f] g + f (\theta H + \varphi X) g 
  \eqno(64)$$
with help of (51). From the commutation relations of the Cartan-Maurer forms 
with the coordinates of the quantum plane, we can compute the corresponding 
relations of $\theta$ and $\varphi$ with functions of the coordinates. From 
(43) we have 
$$(x^m y^n) \theta = q^{-(m+n)} \theta (x^m y^n) + 
  (q^{-n} - 1) \varphi x^m y^{n-1} \qquad 
  (x^m y^n) \varphi = \varphi (x^m y^n). \eqno(65)$$
Inserting (65) in (64) and equating coefficients of the Cartan-Maurer forms, 
we get, for example, 
$$H(f g) = (H f) g + q^{-N} f (H g). \eqno(66) $$
Consequently, we have the coproduct 
$$\Delta(H) = H \otimes 1 + q^{- N} \otimes H $$
$$\Delta(X) = X \otimes 1 + 1 \otimes X + (q^{-1} - 1) X \otimes H. \eqno(67)$$
The counit and coinverse may be calculated by using the axioms of Hopf 
algebra: 
$$m(\epsilon \otimes \mbox{id}) \Delta(u) = u \qquad 
  m(\mbox{id} \otimes \kappa) \Delta(u) = \epsilon(u). \eqno(68)$$
So we have 
$$\epsilon(H) = 0 = \epsilon(X) \eqno(69)$$
$$\kappa(H) = - q^N H \qquad \kappa(X) = - X q^N. \eqno(70)$$

\vspace*{0.3cm}\noindent
{\bf 6. The dual of the Hopf algebra $\cal A$}

\noindent
In this section, in order to obtain the dual of the Hopf algebra ${\cal A}$ 
defined in section 2, we shall use the method of Ref. 12.

A pairing between two vector spaces ${\cal U}$ and ${\cal A}$ is a bilinear 
mapping $<,> : {\cal U} \mbox{x} {\cal A} \longrightarrow {\cal C}, \quad 
  (u,a) \mapsto <u,a>$. We say that the pairing is non-degenerate if 
$$<u,a> = 0 (\forall a \in {\cal A}) ~\Longrightarrow~ u = 0$$
and 
$$<u,a> = 0 (\forall u \in {\cal U}) ~\Longrightarrow~ a = 0.$$
Such a pairing can be extended to a pairing of ${\cal U} \otimes {\cal U}$ 
and ${\cal A} \otimes {\cal A}$ by 
$$<u \otimes v, a \otimes b> = <u,a> <v,b>.$$

Given bialgebras ${\cal U}$ and ${\cal A}$ and a non-degenerate pairing 
$$<,> : {\cal U} \mbox{x} {\cal A} \longrightarrow {\cal C}  \qquad 
  (u,a) \mapsto <u,a> \quad 
  \forall u \in {\cal U} \quad \forall a \in {\cal A} \eqno(71)$$
we say that the bilinear form realizes a duality between ${\cal U}$ and 
${\cal A}$, or that the bialgebras ${\cal U}$ and ${\cal A}$ are in duality, 
if we have 
$$<uv, a> = <u \otimes v, \Delta_{\cal A}> $$
$$<u, ab> = <\Delta_{\cal U}(u), a \otimes b> $$
$$<1_{\cal U}, a> = \epsilon_{\cal A}(a) \eqno(72)$$
and 
$$<u, 1_{\cal A}> = \epsilon_{\cal U}(u) $$
for all $u, v \in {\cal U}$ and $a, b \in {\cal A}$. 

If, in addition, ${\cal U}$ and ${\cal A}$ are Hopf algebras with coinverse 
$\kappa$, then they are said to be in duality if the underlying bialgebras are 
in duality and if, moreover, we have 
$$<\kappa_{\cal U}(u), a> = <u, \kappa_{\cal A}(a)> \qquad 
\forall u \in {\cal U} \quad a \in {\cal A}. \eqno(73)$$

It is enough to define the pairing (71) between the generating elements of the 
two algebras. Pairing for any other elements of ${\cal U}$ and ${\cal A}$ 
follows from relations (72) and the bilinear form inherited by the tensor 
product. For example, for 
$$\Delta_{\cal U}(u) = \sum_k u_k' \otimes u_k'',$$
we have 
$$<u, ab> = <\Delta_{\cal U}(u), a \otimes b> = 
  \sum_k <u_k', a> <u_k'', b> $$

As a Hopf algebra ${\cal A}$ is generated by the elements $x, y$ and a basis 
is given by all monomials of the form 
$$f = x^m y^n $$
where $m, n \in {\cal Z}_+$. Let us denote the dual algebra by 
${\cal U}_q$ and its generating elements by $A$ and $B$. 

The pairing is defined through the tangent vectors as follows 
$$<A, f> = m \delta_{n,0} $$
$$<B, f> = \delta_{n,1}. \eqno(74)$$
We also have 
$$<1_{\cal U}, f> = \epsilon_{\cal A}(f) = \delta_{n,0}. \eqno(75)$$
Using the defining relations one gets 
$$<AB, f> = (m + 1) \delta_{n,1} \eqno(75\mbox{a})$$
and 
$$<BA, f> = m \delta_{n,1} \eqno(75\mbox{b})$$
where differentiation is from the right as this is most suitable for 
differentiation in this basis. Thus one obtains the commutation relation in 
the algebra ${\cal U}_q$ dual to ${\cal A}$ as: 
$$AB = BA + B. \eqno(76)$$
The Hopf algebra structure of this algebra can be deduced by using the 
duality. The coproduct of the elements of the dual algebra is given by 
$$\Delta_{\cal U}(A) = A \otimes 1_{\cal U} + 1_{\cal U} \otimes A $$
$$\Delta_{\cal U}(B) = B \otimes q^A + 1_{\cal U} \otimes B. \eqno(77)$$
The counity is given by 
$$\epsilon_{\cal U}(A) = 0 \qquad \epsilon_{\cal U}(B) = 0. \eqno(78)$$
The coinverse is given as 
$$\kappa_{\cal U}(A) = - A \qquad \kappa_{\cal U}(B) = - B q^{-A}. \eqno(79)$$

We can now transform this algebra to the form obtained in section 5 by making 
the following definitions: 
$$H = {{1_{\cal U} - q^A}\over {1 - q^{-1}}} \qquad X = B \eqno(80)$$
which are consistent with the commutation relation and the Hopf structures. 

\vfill\eject
\noindent
{\bf 7. Conclusion}

\noindent
To conclude, we introduce here the commutation relations between the 
coordinates of the quantum plane and their partial derivatives and thus 
illustrate the connection between the relations in section 5, and the 
relations which will be now obtained. 

To proceed, let us obtain the relations of the coordinates with their partial 
derivatives. We know that the exterior differential {\sf d} can be expressed 
in the form 
$${\sf d} f = ({\sf d} x ~\partial_x + {\sf d} y ~\partial_y) f. 
  \eqno(81)$$
Then, for example, 
\begin{eqnarray*}
{\sf d} (x f) 
& = & {\sf d} x ~f + x ~{\sf d} f\\
& = & {\sf d} x ~(1 + q^{-1} x \partial_x) f + {\sf d} y ~x \partial_y f \\
& = & ({\sf d} x ~\partial_x x + {\sf d} y ~\partial_y x) f
\end{eqnarray*}
so that 
$$\partial_x x = 1 + q^{-1} x \partial_x 
  \qquad \partial_x y = q^{-1} y \partial_x $$
$$\partial_y x = x \partial_y \qquad 
  \partial_y y = 1 + q^{-1} y \partial_y + (q^{-1} - 1) x \partial_x. 
  \eqno(82)$$
The Hopf algebra structure for $\partial$ is given by 
$$\Delta(\partial_x) = \partial_x \otimes \partial_x \qquad 
  \Delta(\partial_y) = \partial_y \otimes 1 + \partial_x \otimes \partial_y $$
$$\epsilon(\partial_x) = 1, \qquad \epsilon(\partial_y) = 0 \eqno(83)$$ 
$$\kappa(\partial_x) = \partial_x^{-1} \qquad 
  \kappa(\partial_y) = - \partial_x^{-1} \partial_y $$
provided that the formal inverse $\partial_x^{-1}$ exists. 

We know, from section 5, that the exterior differential {\sf d} can be 
expressed in the form (51), which we repeat here, 
$${\sf d} f = (\theta H + \varphi X) f. \eqno(84)$$
Considering (81) together (84) and using (50) one has 
$$H \equiv x \partial_x + y \partial_y \qquad 
  X \equiv \partial_y. \eqno(85)$$
Using the relations (82) and (32) one can check that the relation of the 
generators in (85) coincide with (55). It can also be verified that, the 
action of the generators in (85) on the coordinates coincide with (58). 

We finally introduce complex notation with a single variable 
$z = x + {\bf i} y$ where $x$ and $y$ are the generators of the $q$-plane and 
${\bf i}^2 = - 1$. Then the elements 
$$z \qquad \bar{z} = x - {\bf i}y \qquad  
  {\sf d} z = {\sf d} x + {\bf i} {\sf d} y \qquad 
  {\sf d} {\bar z} = {\sf d} x - {\bf i} {\sf d} y \eqno(86)$$
form the basis in the algebra $\Gamma$. These elements obey the following 
commutation relations: 
$$z ~{\sf d} z = q^{-1} {\sf d} z ~z \qquad 
  (\bar{z} ~{\sf d} \bar{z} = q^{-1} {\sf d} \bar{z} ~\bar{z})$$
$$z ~{\sf d}^2 z = q^{-1} {\sf d}^2 z ~z \qquad 
  (\bar{z} ~{\sf d}^2 \bar{z} = q^{-1} {\sf d}^2 \bar{z} ~\bar{z})$$
$${\sf d} z ~{\sf d}^2 z = q^{-2} {\sf d}^2 z ~{\sf d} z \qquad 
({\sf d}\bar{z} ~{\sf d}^2 \bar{z} = q^{-2}{\sf d}^2 \bar{z} ~{\sf d}\bar{z})$$
$$({\sf d} z)^3 = 0 = ({\sf d} \bar{z})^3. \eqno(87)$$
Note that these relations are the same with those of Ref. 11 except that in 
our case $z^3$ need not be zero.

The Z$_3$-graded noncommutative differential geometry we have constructed 
satisfies all expectations for such a structure. In particular all Hopf 
algebra axioms are satisfied without any modification. 

\noindent
{\bf Acknowledgment}

\noindent
This work was supported in part by {\bf T. B. T. A. K.} the 
Turkish Scientific and Technical Research Council. 

\vfill\eject

\end{document}